\documentclass[12pt]{article}

\usepackage{latexsym}
\usepackage{amsfonts}
\usepackage{amssymb}
\usepackage{amsmath}
\usepackage{mathrsfs}
\usepackage{bbm}
\usepackage{dsfont}
\usepackage{tikz-cd}

\usepackage{color}
\definecolor{indigo}{RGB}{0,0,120}
\usepackage[colorlinks=true,linkcolor=indigo,citecolor=blue,urlcolor=indigo]{hyperref}

\newcommand{\be}{\begin{equation}}
\newcommand{\ee}{\end{equation}}
\newcommand{\bea}{\begin{eqnarray}}
\newcommand{\eea}{\end{eqnarray}}
\newcommand{\bean}{\begin{eqnarray*}}
\newcommand{\eean}{\end{eqnarray*}}
\newcommand{\brray}{\begin{array}}
\newcommand{\erray}{\end{array}}

\newtheorem{dfn}{Definition}[section]
\newtheorem{thm}[dfn]{Theorem}
\newtheorem{lmma}[dfn]{Lemma}
\newtheorem{ppsn}[dfn]{Proposition}
\newtheorem{crlre}[dfn]{Corollary}
\newtheorem{xmpl}[dfn]{Example}
\newtheorem{rmrk}[dfn]{Remark}

\newcommand{\bdfn}{\begin{dfn}\rm}
\newcommand{\bthm}{\begin{thm}}
\newcommand{\blmma}{\begin{lmma}}
\newcommand{\bppsn}{\begin{ppsn}}
\newcommand{\bcrlre}{\begin{crlre}}
\newcommand{\bxmpl}{\begin{xmpl}}
\newcommand{\brmrk}{\begin{rmrk}\rm}

\newcommand{\edfn}{\end{dfn}}
\newcommand{\ethm}{\end{thm}}
\newcommand{\elmma}{\end{lmma}}
\newcommand{\eppsn}{\end{ppsn}}
\newcommand{\ecrlre}{\end{crlre}}
\newcommand{\exmpl}{\end{xmpl}}
\newcommand{\ermrk}{\end{rmrk}}

\title{Opposite product system for the multiparameter CAR flows}
\author{Anbu Arjunan}

\begin{document}
\maketitle
\begin{abstract}
  We consider the multiparameter CAR flows and describe its opposite. We also characterize the symmeticity of CAR flows in terms of associated isometric representations.
 \end{abstract}
\noindent {\bf AMS Classification No. :} {Primary 46L55; Secondary 46L99.}  
\section{Introduction}
Let $P$ be a closed convex cone in $\mathbb{R}^d$. We assume that $P-P=\mathbb{R}^d$ and $P\cap -P=\{0\}$. Let $V$ be a pure isometric representation of $P$ and let $\alpha$ be the CCR flow associated to the isometric representation $V$. 
The author in \cite{R19} have shown that the CCR flow is not cocycle conjugate to the CAR flow when the isometric representation $V$ is proper. 
The product system associated with the CAR flow is not decomposable in general; see \cite{arjunan2020decomposability}.
It was shown in \cite{sundar2020asymmetric} that $\alpha$ is cocycle conjugate to $ \alpha^{\text{op}}$ if and only if $V$ is unitary equivalent to its opposite $V^{\text{op}}$.
This result uses the characterization of decomposable product system which admits a unit; see \cite{sundar2019arvesons}.
It is natural to ask whether the analogous result holds true for the multiparameter CAR flows. In this article we answer this question affirmatively; see Theorem \ref{main}.
 We will achieve this by identifying the opposite of the product system for a CAR flow with the product system for an appropriate CAR flow. Also we will also use this to study the symmetricity of the CAR flows.
\section{Preliminaries}
 Let $H$ be a Hilbert space and let $\Gamma_a(H)$ be the antisymmetric Fock space over $H$.
 For $f\in H$, define a bounded operator $a(f)^*$ on $\Gamma_a(H)$ as 
 \begin{align*}
  a(f)^*(\Omega)&=f\text{ and }\\
  a(f)^*(h_1\wedge h_2\wedge...\wedge h_n)&=f\wedge h_1\wedge h_2\wedge...\wedge h_n
 \end{align*}
 where $\Omega$ is the vacuum vector of $\Gamma_a(H)$ and $h_1\wedge h_2\wedge...\wedge h_n$ is an arbitrary antisymmetric elementary tensor element with $h_1, h_2,..., h_n\in H$ and $n\geq 1$. Let $a(f)$ be the adjoint of $a(f)^*$. The operators $a(f)^*$ and $a(f)$  are called the creation and the annihilation operator associated to a vector $f$. 
 
 By an isometric representation of $P$ on a Hilbert space $H$, we mean a strongly continuous map $V:P\to B(H)$ such that each $V_x$ is an isometry and $V_xV_y=V_{x+y}$ for each $x,y\in P$. 
 For a given isometric representation $V:P\to B(H)$, there exists a unique $E_0$-semigroup, denoted by $\beta^V$, on $\Gamma_a(H)$ satisfying 
 \[\beta^V_x(a(f))=a(V_xf)\text{ for each }f\in H.\]
 This $E_0$-semigroup $\beta^V$ is called the CAR flow associated to the isometric representation $V$; see \cite{R19}.
 
 Let $H$ and $K$ be Hilbert spaces. For an isometry $W:H\to K$, there exists a unique bounded operator $\Gamma_a(W)$, called the second quantization of $W$, from $\Gamma_a(H)$ to $\Gamma_a(K)$, satisfying 
\begin{align*}
\Gamma_a(W)(\Omega)&=\Omega, \text{ and }\\
\Gamma_a(W)(f_1\wedge f_2\wedge...\wedge f_n)&=Wf_1\wedge Wf_2\wedge...\wedge Wf_n ,
 \end{align*}
 where $\Omega$ is the vacuum vector in the appropriate antisymmetric Fock space and $f_1\wedge f_2\wedge...\wedge f_n$ is any antisymmetric elementary tensor element with $f_1, f_2,..., f_n\in H$ and $n\geq 1$. 
\section{Opposite product sysem for a CAR flow} 
Let $V$ be a pure isometric representation of $P$ on a Hilbert space $H$. Let $\beta^V$ be the CAR flow associated to the isometric representation $V$ and denote its concrete product system by $\mathcal{E}_{\beta^V}$.  
Set $E^V(x)=\Gamma_a(\text{Ker}(V_x^*))$.
Consider the set $E^V$ as 
\[E^V=\{(x,f):x\in \Omega\text{ and }f\in E^V(x)\}.\]
Since $E^V$ is a Borel subset of $\Omega\times \Gamma_a(H)$, $E^V$ is a standard Borel space. Define a multiplication $.$ on $E^V$ as 
\[(x,f).(y,f):=(x+y, f\otimes \Gamma_a(V_x)g)\]
for every $(x,f),(y,f)\in E^V$. $E_V$ equipped with the above multiplication defines a product system structure over $\Omega$. We define another multiplication $\circ$ on $E^V$ as 
\[(x,f)\circ(y,f):=(x+y, g\otimes \Gamma_a(V_y)f).\]
Then the pair $(E^V,\circ)$ also has a structure of product system over $\Omega$, called the opposite product system for $(E^V,.)$, denoted by $(E^V)^{\text{op}}$.

Let $x\in \Omega$ and let $f\in E^V(x)$ be given. Define a bounded operator $T_f$ on $\Gamma_a(H)$ as 
\[T_f\eta=f\otimes \Gamma_a(V_x)\eta,\text{ for every }\eta\in \Gamma_a(H).\]
Then we have the following lemma.
\begin{lmma}\label{keylemma} The map  $\theta:E^V\ni (x,f)\mapsto (x,T_{f})\in \mathcal{E}_{\beta^V}$ is an isomorphism as product systems.
\end{lmma}
\begin{prf}
 Let $(x,f),(y,g)\in E^V$ be given. Since $T_fT_g=T_{f\otimes \Gamma_a(V_x)g}$, it follows that $\theta(x,f)\theta(y,g)=\theta((x,f)(y,g))$. For each $x\in \Omega$, the restriction of $\theta$ to $E^V(x)$, $\theta|_{E^V(x)}:E^V(x)\to \mathcal{E}_{\beta^V}(x)$ is a unitary. For let $f,g\in E^V(x)$ be given. Note that $T_g^*T_f=\langle f,g\rangle 1_{E^V(x)}$ and $T_f\in \mathcal{E}_{\beta^V}(x)$.
 This implies that the map $E^V(x)\ni f\mapsto T_f\in \mathcal{E}_{\beta^V}(x)$ is an isometry. To prove that the map is a unitary it suffices to show that whenever $T\in \mathcal{E}_{\beta^V}(x)$ such that $\langle T_f,T\rangle=0$ for all $f$, then $T=0$.
 Since the linear span of the set $\{f\otimes \Gamma_a(V_x)\eta: f\in E^V(x)\text{ and }\eta\in \Gamma_a(H)\}$ is dense in $\Gamma_a(H)$, we see that $T=0$.
 
 Since $E^V$ and $\mathcal{E}_{\beta^V}$ are standard Borel spaces and the restriction of $\theta$ to each fibre is a unitary, it follows that the map $\theta$ is a Borel isomorphism and hence it is a isomorphism as product systems by \cite{ArvCA}.
 \hfill $\Box$
\end{prf}

Let us recall the opposite isometric representation $V^{\text{op}}$ for the given isometric representation $V$ considered in \cite{sundar2019arvesons}.
Let $U$ be a minimal unitary dilation of $V$. More precisely, there exists a Hilbert space $\widetilde{H}$ containing $H$ as a subspace  and a unitary representation $U$ of $\mathbb{R}^d$ on a Hilbert $\widetilde{H}$ such that the following conditions hold.
\begin{enumerate}
 \item For $x\in P$, $U_x\xi=V_x\xi$.
 \item The set $\cup_{x\in P}U_x^*H$ is dense in $\widetilde{H}$.
\end{enumerate}
Note that for $x\in P$, $K=H^{\perp}$ is invariant under $U_x$. For $x\in P$, define $V^{\text{op}}_x$ on $K$ to be the restriction of $U_{-x}$ to $K$ i.e. $V^{\text{op}}_x:=U_{-x}|_{K}$. Then $V^{\text{op}}:=\{V^{\text{op}}_x\}_{x\in P}$ is an isometric representation of $P$, called the opposite isometric representation for $V$. This isometric representation $V^{\text{op}}$ is pure \cite[Proposition 3.2]{sundar2020asymmetric}.
\begin{ppsn}\label{keyprop}
 The map $\phi:(E^V)^{\text{op}}\ni (x,f)\mapsto (x,\Gamma_a(U_{-x})f)\in E^{V^{\text{op}}}$ is an isomorphism as product systems.
 \end{ppsn}
 \begin{prf}
 For each $x\in \Omega$, the map $\text{Ker}(V_x^*)h\mapsto U_{-x}h\in \text{Ker}((V_x^{\text{op}})^*)$ is a unitary; see the proof of \cite[Proposition 3.2]{sundar2020asymmetric}. Then it follows that the map $\phi:(E^V)^{\text{op}}\ni (x,f)\mapsto (x,\Gamma_a(U_{-x})f)\in E^{V^{\text{op}}}$ is a continuous bijection and its inverse is given by $E^{V^{\text{op}}}\ni (x,\xi)\mapsto (x,\Gamma_a(U_{x})\xi)\in(E^V)^{\text{op}}$. Hence it is a Borel isomorphism by \cite{ArvCA}. Now it remains to show that $\phi$ follows product system structure.
  Let $(x,f),(y,g)\in E^V$ be given. Then we have
  \begin{align*}
   \phi((x,f)(y,g))&=\phi(x+y,f\otimes \Gamma_a(V_x)g)\\
   &=(x+y,\Gamma_a(U_{-(x+y)})(f\otimes \Gamma_a(V_x)g))\\
   &=(x+y,\Gamma_a(U_{-(x+y)}) \Gamma_a(V_x)g\otimes \Gamma_a(U_{-(x+y)})f)\\
   &=(x+y,\Gamma_a(U_{-y})g\otimes \Gamma_a(U_{-y})\Gamma_a(U_{-x})f)\\
   &=(x+y,\Gamma_a(U_{-y})g\otimes \Gamma_a(V_y^{\text{op}})\Gamma_a(U_{-x})f)\\
   &=(y,\Gamma_a(U_{-y})g)(x,\Gamma_a(U_{-x})f)\\
   &=\phi(y,g)\phi(x,f).
  \end{align*}
Hence the map $\phi$ is an isomorphism as product systems.
\hfill $\Box$
 \end{prf}
 
 Let $\mathcal{E}_{\beta^V}$ be the concrete product system for $\beta^V$ and let $\mathcal{E}_{\beta^V}^{\text{op}}$ be its opposite product system. By \cite[Theorem 3.14]{MS}, there exists an $E_0$-semigroup denoted by $(\beta^V)^{\text{op}}$ such that $\mathcal{E}_{\beta^V}^{\text{op}}$ is isomorphic to $\mathcal{E}_{(\beta^V)^{\text{op}}}$.
\begin{crlre}\label{oppositecar}
 An $E_0$-semigroup $(\beta^V)^{\text{op}}$ is cocycle conjugate to $\beta^{V^{\text{op}}}$.
  \end{crlre}
  \begin{prf}
   By Proposition \ref{keyprop} and Lemma \ref{keylemma}, we conclude that $(E^V)^{\text{op}}$ is isomorphic to $\mathcal{E}_{\beta^{V^{\text{op}}}}$. This implies that the product system $\mathcal{E}_{(\beta^V)^{\text{op}}}$ is isomorphic to $\mathcal{E}_{\beta^{V^{\text{op}}}}$ by Lemma \ref{keylemma}.
   Then by \cite[Theorem 2.9]{MS}, we have $(\beta^V)^{\text{op}}$ is cocycle conjugate to $\beta^{V^{\text{op}}}$. 
   \hfill $\Box$
  \end{prf}
\begin{rmrk}
   The above corollary implies that the opposite of a CAR flow over $P$ is again a CAR flow over $P$.
  \end{rmrk}
\begin{thm} \label{main} Let $\beta^V$ be the CAR flow associated to an isometric representation $V$. Then the following are equivalent.
 \begin{enumerate}
  \item The CAR flow $\beta^V$ is cocycle conjugate to its opposite $(\beta^V)^{\text{op}}$
  \item The isometric representation $V$ is unitary equivalent to its opposite $V^{\text{op}}$.
 \end{enumerate}
\end{thm}
\begin{prf}
Proof follows from \cite[Proposition 4.7]{R19} and Corollary \ref{oppositecar}.
\hfill $\Box$
\end{prf}
\section{Examples for symmetric and asymmetric CAR flows}
By a $P$-module we mean a non-empty closed subset $A$ of $\mathbb{R}^d$ such that $A+P\subseteq A$.
Let $A$ be a $P$-module. For $x\in P$, define an operator $V^A_x$ on $L^2(A)$ as
\[(V_x^{A}f)(y)=\begin{cases}
 f(y-x)  & \mbox{ if } y-x\in A,\cr
   
    0 &  \mbox{ if }  y-x\notin A.
         \end{cases}\]
for each $f\in L^2(A)$. Then the family $\{V_x^A\}_{x\in P}$ is an isometric representation of $P$. 
\begin{ppsn}\label{sundarresult} (See \cite[Proposition 3.4]{sundar2020asymmetric})
We have the following.
\begin{enumerate}
 \item The isometric representation $(V^A)^{\text{op}}$ is unitary equivalent to $V^A$.
 \item There exists an element $z\in \mathbb{R}^d$ such that $A=-(\text{int}(A)^c)+z$.
\end{enumerate}
 Here $\text{int}(A)$ is the interior of $A$ and $\text{int}(A)^c$ is the complement of $\text{int}(A)$ in $\mathbb{R}^d$.
 \end{ppsn}
Let $\beta^A$ be the CAR flow associated to the isometric representation $V^A$. It follows from Theorem \ref{main} and Proposition \ref{sundarresult} that the CAR flow $\beta^A$ is cocycle conjugate to its opposite $(\beta^A)^{\text{op}}$ if and only if $A=-(\text{int}(A)^c)+z$ for some $z\in\mathbb{R}^d$.
\begin{rmrk}
 By considering the existence of such $P$-modules, we can see that there are uncountably many symmetric CAR flows as well as asymmetric CAR flows over $P$.
\end{rmrk}
\section*{Acknowledgment}
 The author would like to thank The Institute of Mathematical Sciences for the Institute Postdoctoral fellowship.
\bibliography{reference}
 \bibliographystyle{amsplain}

\end{document}